\documentclass[journal]{IEEEtran}

\usepackage{comment}

\usepackage{cite}

\ifCLASSINFOpdf
  \usepackage[pdftex]{graphicx}
  \graphicspath{{../../Figures/}{../jpeg/}}
\else
  \usepackage[dvips]{graphicx}
  \graphicspath{{../../Figures/}{../jpeg/}}
  \DeclareGraphicsExtensions{.eps}
\fi
\usepackage[cmex10]{amsmath}
\ifCLASSOPTIONcompsoc
  \usepackage[caption=false,font=normalsize,labelfont=sf,textfont=sf]{subfig}
\else
  \usepackage[caption=false,font=footnotesize]{subfig}
\fi
\usepackage[hyphens]{url}

\usepackage{framed}
\usepackage{amssymb}
\usepackage{color}
\usepackage[T1]{fontenc}

\begin{document}
%
\title{Data Centers as Dispatchable Loads to \\ Harness Stranded Power}
%
%
%


\author{Kibaek~Kim,~\IEEEmembership{Member,~IEEE,}
Fan~Yang,~\IEEEmembership{Member,~IEEE,}
Victor~M.~Zavala,~\IEEEmembership{Member,~IEEE,}
and~Andrew~A.~Chien,~\IEEEmembership{Fellow,~IEEE}
\thanks{Kibaek Kim is with the Mathematics and Computer Science (MCS) Division, Argonne National Laboratory, Lemont, IL 60439, Email: kimk@anl.gov.}
\thanks{Fan Yang and Andrew A. Chien are with the Department of Computer Science at the University of Chicago.}
\thanks{Victor M. Zavala is with the Department of Chemical and Biological Engineering at the University of Wisconsin-Madison.}
\thanks{Andrew A. Chien and Victor M. Zavala are also affiliated with the MCS Division at Argonne National Laboratory.}
}
\maketitle

\begin{abstract}
We analyze how both traditional data center integration and
dispatchable load integration affect power grid efficiency.
We use detailed network models, parallel optimization solvers, and thousands of renewable
generation scenarios to perform our analysis.  
Our analysis reveals that significant spillage and stranded power
will be observed in power grids as wind power levels are increased.  A counter-intuitive finding is that collocating data centers with inflexible loads next to 
wind farms has limited impacts on renewable portfolio standard (RPS) goals because it 
provides limited system-level flexibility and can in fact increase stranded power and fossil-fueled generation.  In contrast,
optimally placing data centers that are dispatchable (with flexible loads) provides 
system-wide flexibility, reduces stranded power, and improves efficiency.  In short, optimally placed dispatchable computing loads can enable better scaling to
high RPS.  We show that these dispatchable computing loads are powered
to 60$\sim$80\% of their requested capacity, indicating that there are significant economic incentives provided by stranded power.
\end{abstract}

\begin{IEEEkeywords}
renewable power, green computing, power grid, energy markets, renewable portfolio
standard, cloud computing
\end{IEEEkeywords}

%
\IEEEpeerreviewmaketitle




\section{Introduction}
\label{section:introduction}

Over the past two decades, a growing consensus on climate change due to
anthropogenic carbon \cite{IPCC-CC2014,Gore06} has emerged.  In response,
efforts have expanded worldwide to reduce the amount of carbon being released
into the atmosphere \cite{Kyoto1997,Paris2015}.  Information and computing
technologies (ICT) emissions, estimated at 2\% in 2010 of global emissions
\cite{ICT-Footprint10}, are among the most rapidly growing.  In fact, recent
estimates put ICT emissions at 8\% of electric power in 2016, growing to 13\%
by 2027 \cite{ICT-Footprint10,GreenPeaceCC15}.  The rise of cloud computing has
led to growing concerns about carbon emissions from data centers
\cite{GP-cclean14,GreenPeaceCC15}, spawning research on data-center energy
efficiency \cite{Barroso09,Patel05} and how to exploit renewable power
\cite{GoogleGreen,AppleNC,goiri2013parasol,berral2014building}.  A recent
strategy pursued by several ``hyperscaler'' internet companies is to purchase
wind-power offsets as part of ``long-term purchase'' contracts
\cite{ITwindcontract}.  Another well-studied research topic is the optimization
of data-center site selection based on cost and on exploitation of renewable
power \cite{wang15grid,le2011intelligent}.  To the best of our knowledge, all
such studies consider costs and benefits from a cloud computing operator
perspective, seeking to maximize data center revenue, reduce total cost of
ownership (TCO), and maintain high data-center availability.  In this paper, we
consider the impact of the addition of new data centers and collocated
renewables on system-level resilience, efficiency, and flexibility of the power
grid. 

Closely related work includes study of data-center demand-response (DCDR)
\cite{Wierman-DCDR14} and how to implement it in data centers by affecting
scheduling and providing economic incentives to data center tenants
\cite{Wierman-Greening15,Wierman-Pricing14}. These models employ grid economic
incentives and aim to reduce data center load by 15$\sim $20\% on demand. A
fundamental challenge for such approaches is that data-center operators are not
inclined to participate - even with significant pricing incentives, because
power cost is typically $<$10\% of data-center TCO, and demand-response
requires significant new complexity.  In contrast, we consider a new model
where the data center is a ``dispatchable load,'' based on new computer science
approaches to create flexible computing loads. This produces two advantages:
(1) much larger dynamic load range (100\% of the data-center use vs. 15\%) and
(2) elimination of the need to convince traditional data center operators to
undertake DCDR.  \footnote{This framing does not preclude the possibilities
that these dispatachable data centers will receive favorable economic treatment
from the grid, since they can provide a significant positive benefit.}

Ambitious ``renewable portfolio standards'' (RPS) goals for renewable power as
a fraction of overall power have been widely adopted.  Midwest examples from
the Mid-continent Independent System Operator (MISO) system include Illinois
(25\% by 2025) to Minnesota (25$\sim$31\% by 2025).  California and New York
have adopted a 50\% goal by 2030 \cite{rps-california,Megerian15}.  Obama's
``Clean Power Plan,'' issued August 2015, calls for a national 32\% reduction
in electric power carbon emissions by 2030, with renewable power as a critical
element.  And, the U.S. Department of Energy landmark report ``Wind Vision
2015'' describes how the United States can achieve a 35\% RPS for wind alone by
2050, a big jump from a combined solar and wind RPS of 5.2\% in 2014
\cite{renewable2014}.  In this plan, some regions such as the Midwest and Texas
achieve RPS over 50\% by 2050.  These ambitious and transformative goals pose
serious power grid challenges, including the ability to achieve ``merit
order,'' efficiency, stability, and resiliency. 

With dual goals of high RPS in the power grid and support for large-scale
computing, we address three questions:

\begin{enumerate}
  \item What is the impact of data centers on the future power grid?
  \item Should renewable generators be collocated with data centers?
  \item Can we enable growing cloud computing and renewable penetration? 
\end{enumerate}

To explore these questions, we developed a computational framework that
combines a detailed power grid network models and cutting-edge parallel
optimization solvers to identify optimal designs that remain resilient in the
face of a wide range of operational scenarios. We consider different settings
that are variations of today's Western Electricity Coordinating Council (WECC)
grid: a base system, the addition of twenty 200 MW data centers at random
locations, the same twenty data centers with collocated renewables and
inflexible loads, and the same data centers with dispatchable loads (volatile
cloud computing) and optimally placed in the system.  Such dispatchable loads
can enable grid efficiency at higher RPS levels because they enable system-wide
flexibility.  We consider optimizing the location of these dispatchable loads
and the resulting impact on power grid efficiency.  For each case we
characterize the impact of increasing RPS levels and data centers on grid
system cost, stranded power, stability, and the capacity achieved by the data
centers. Our findings include the following.  
\begin{itemize}

  \item Significant spillage and stranded power exists in current power grids
and increase with higher RPS levels.

  \item Collocating wind farms and data centers can in fact be harmful to RPS
goals, increasing stranded power and thermal generation.

  \item Dispatchable computing loads reduce stranded power and enable higher
RPS.  

  \item Optimizing dispatchable load locations decreases system-wide social
and powers the dispatchable data centers to 60$\sim$80\% of capacity. 

\end{itemize}

The rest of the paper is organized as follows.  In Section \ref{section:loads}
we introduce concepts of stranded power and dispatchable computing loads.
Section \ref{section:models} outlines the optimization models, followed by
experiments in Section \ref{section:experiments}. In
Section~\ref{section:summary} we summarize our results and briefly discuss
future work.


\section{Stranded Power and Dispatchable Computing}
\label{section:loads}

\subsection{Stranded Power}\label{sec:stranded}

Power system operators must balance power flows across the network.  Generators
offer their generation capability to the grid in real time (every 5 to 12
minutes), and the grid dispatches generation based on the demand and
transmission.  Intrinsic variability of renewable generation (wind, solar,
etc.) creates major challenges for power dispatch.  Despite best efforts to
match generation and demand, in the process of ensuring reliable power,
oversupply and transmission congestion can prevent generated power from
reaching loads at certain times.  Power grids call this excess power
\emph{spillage}, ``curtailment,'' or ``down dispatching.''

\begin{figure}[htb]
  \centering
  \includegraphics[width=0.45\textwidth]{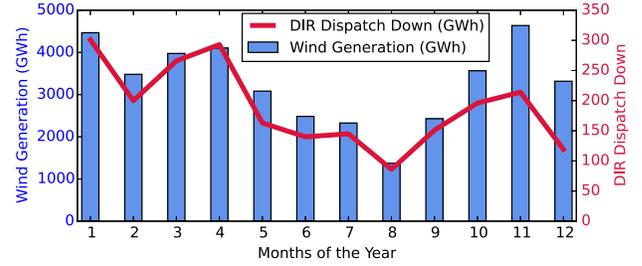}
  \caption{Wind generation and down-dispatching (spillage) of MISO in 2014.}
  \label{fig:miso_dispatch}
\end{figure}

Figure \ref{fig:miso_dispatch} shows the monthly wind generation and
down-dispatched wind power (spillage) of the MISO system.  Almost 7\% of wind
generation is curtailed because of transmission congestion.  The total
down-dispatched power of MISO for 2014 was about 2.2 terawatt-hours,
corresponding to a 183 MW sustained rate.  Comparable waste also exists in
other independent system operators (ISO), including the Eastern Region
Coordinating District of Texas, California ISO (CAISO) \cite{E3-Fifty-Percent},
and many European countries such as Denmark, Germany, Ireland, and Italy
\cite{Lew2013}.  The amount of waste is projected to increase with higher RPS
levels \cite{Megerian15}.

Modern energy markets dispatch generation by assigning locational marginal
prices (LMPs) that vary across generation sites, grid nodes, and time
intervals.  In situations of oversupply or transmission congestion, power
prices can be low or even negative, causing power generators to dump power
(spillage) or deliver it and pay the grid to take it.\footnote{Delivering power
into the grid can be tied to ``production tax credits,'' a financial incentive
to sell power to the grid even at negative price.} Consequently, spillage can
be significantly less than total uneconomic generation.

\begin{figure}[htb]
  \centering
  \includegraphics[width=0.45\textwidth]{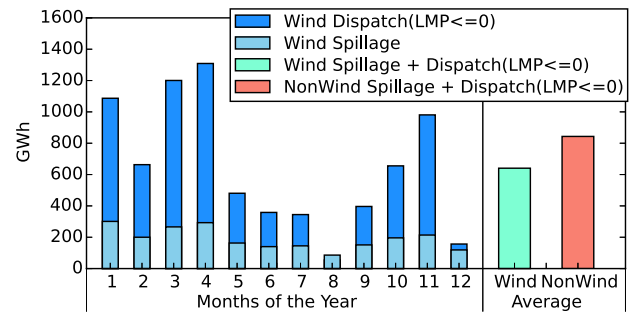}
  \caption{MISO stranded power (wind) vs. stranded power (nonwind) in 2014.}
  \label{fig:miso_stranded_power}
\end{figure}

We define \emph{stranded power} as all offered generation with no economic
value, thus including both spillage and delivered power with zero or negative
LMP.  Figure \ref{fig:miso_stranded_power} presents MISO's stranded power in
2014, breaking it down by month and type.  It also compares the average
stranded power from wind and nonwind generation.  Overall stranded wind power,
the sum of wind spillage and noneconomic wind dispatch (LMP$\leq$0), is about
7.7 TWh for 2014.  Interestingly, nonwind sites, mostly thermal generators,
have 10.1 TWh of stranded power. Because 90\% of grid power is thermal, however, 
the thermal stranded power is approximately 8 times less by percentage.

\subsection{Dispatchable Loads}

We define \emph{dispatchable loads} as adjustable (flexible) demands that are
dispatched in real time by the power grid.  Adjustable at each dispatch
interval, dispatchable loads are an ideal way of demand response and appear
promise to reduce congestion and stranded power due to generation variability.
Key properties of dispatchable loads include the following.
\begin{enumerate}
  \item Grid control can increase consumption to limit.
  \item Grid control can decrease consumption to zero.
  \item Grid control can be exercised at the dispatch interval (instantly).
\end{enumerate}

A less obvious benefit of dispatchable loads is their ability to control
network flows (operators currently have limited control over flows as
electricity follows paths of least resistance).  

We consider one possible type of dispatchable load: computing.  Operating on
a short time scale, computing has the potential to be agile yet still
productive. We call such intermittent computing resources zero-carbon clouds
(ZCClouds) and have described and analyzed several possible models
\cite{chien2015zero,yang2016zccloud}.  For example, ZCClouds can be computing
hardware deployed in shipping containers and directly connected to a wind farm
or at a key transmission bottleneck in the power grid.  ZCClouds can transform
stranded power into computing power with short latency (in seconds) and
can be easily turned on or shut down according to stranded power availability.
Possible uses include data-center workloads such as big data analysis, machine
learning, or high-performance computing.  The uptime and capabilities of
intermittent computing resources deployed as dispatchable loads are determined
by the quanity and temporal distribution of stranded power. 

The idea of intermittent (or ``volatile'') computing resources is of growing
interest.  Cloud providers have begun to provide unreliable/revokable resources
including Amazon's spot instances \cite{Spot-Instance-Pricing2013} and Google's
preemptible VM Instance \cite{gce-preemptible}.  Several studies propose
methods to make such resources useful for high-performance computing and more
advanced cloud services \cite{SpotHPDC2015}.  We believe there is a broad
application for intermittent computing resources.

Many other realizations of dispatchable loads are possible, including energy
storage.  A key difference, however, is that dispatchable loads have infinite
capacity (can run forever without filling up) and avoid subsidized economics
(computing services defray their cost by providing services not related to
energy).


\section{Optimization Models}
\label{section:models}

In this section we present two optimization models in order to assess the
benefits of dispatchable loads. We also present various performance metrics to
conduct our analysis.

We begin with the model notation in the following table. The units for power,
energy, and phase angles are megawatts, megawatt-hours, and degrees;
respectively.  

Sets:
\begin{description}[\IEEEusemathlabelsep]
  \item[$\mathcal{D};\mathcal{D}_n$] Demand loads; demand loads at bus $n$
  \item[$\mathcal{G};\mathcal{G}_n$] Generators; generators at bus $n$
  \item[$\mathcal{I};\mathcal{I}_n$] Import points; import points at bus $n$
  \item[$\mathcal{L}$] Transmission lines
  \item[$\mathcal{L}_n^+;\mathcal{L}_n^-$] Transmission lines to bus $n$; lines from bus $n$
  \item[$\mathcal{N}$] Buses
  \item[$\mathcal{R};\mathcal{R}_n$] Renewable generators; Renewable generators at bus $n$
  \item[$\mathcal{T}$] Time periods
  \item[$\mathcal{W}$] Wind-farm locations
  \item[$\Omega$] Scenarios
\end{description}

Parameters:
\begin{description}[\IEEEusemathlabelsep]
  \item[$B_l$] Susceptance of transmission line $l$
  \item[$C_i$] Generation cost of generator $i$ 
  \item[$C_j^d$] Load-shedding penalty at load $j$ 
  \item[$C_i^w$] Spillage penalty at wind farm $i$ 
  \item[$C_i^m$] Spillage penalty at import point $i$ 
  \item[$C_i^r$] Spillage penalty at renewable $i$ 
  \item[$C_n^s$] Value of lost dispatchable load at bus $n$ 
  \item[$D_{jt}$] Demand load of consumer $j$ at time $t$ 
  \item[$F^{max}_l$] Maximum power flow of transmission line $l$ 
  \item[$K$] Maximum number of dispatchable loads
  \item[$M_{it}$] Power production of import $i$ at time $t$ 
  \item[$P^{max}_i$] Maximum power output of generator $i$ 
  \item[$R_{it}$] Power production of renewable $i$ at time $t$ 
  \item[$RU_i$] Ramp-up limit of generator $i$ 
  \item[$RD_i$] Ramp-down limit of generator $i$ 
  \item[$U$] Dispatchable load capacity 
  \item[$W_{wt}$] Power from wind farm $w$ at time $t$ 
  \item[$W_{wt}(\omega)$] Power from wind farm $w$ at time $t$ for scenario $\omega$ 
  \item[$\pi(\omega)$] Probability of wind production scenario $\omega$
  \item[$\Theta_{nt}^{min}$] Minimum phase angle at bus $n$ at time $t$
  \item[$\Theta_{nt}^{max}$] Maximum phase angle at bus $n$ at time $t$
\end{description}

Decision variables:
\begin{description}[\IEEEusemathlabelsep]
  \item[$d_{jt}$] Load shedding at load $j$ at time $t$ 
  \item[$d_{jt}(\omega)$] Load shedding at load $j$ at time $t$ for scenario $\omega$ 
  \item[$f_{lt}$] Power flow of line $l$ at time $t$ 
  \item[$f_{lt}(\omega)$] Power flow of line $l$ at time $t$ for scenario $\omega$ 
  \item[$m_{it}$] Spillage at import $i$ at time $t$ 
  \item[$m_{it}(\omega)$] Spillage at import $i$ at time $t$ for scenario $\omega$ 
  \item[$p_{it}$] Power from generator $i$ at time $t$ 
  \item[$p_{it}(\omega)$] Power from generator $i$ at time $t$ for scenario $\omega$ 
  \item[$r_{it}$] Spillage at renewable $i$ at time $t$ 
  \item[$r_{it}(\omega)$] Spillage at renewable $i$ at time $t$ for scenario $\omega$ 
  \item[$u_{nt}(\omega)$] Dispatchable load served at bus $n$ at time $t$ for scenario $\omega$ 
  \item[$w_{it}$] Spillage at wind farm $i$ at time $t$ 
  \item[$w_{it}(\omega)$] Spillage at wind farm $i$ at time $t$ for scenario $\omega$ 
  \item[$x_{n}$] Number of dispatchable loads installed at bus $n$
  \item[$\theta_{nt}$] Phase angle at bus $n$ at time $t$
  \item[$\theta_{nt}(\omega)$] Phase angle at bus $n$ at time $t$ for scenario $\omega$
\end{description}

\subsection{Economic Dispatch Model}

To assess the economic benefits of dispatchable computing loads, we use the
following economic dispatch (ED) model:

\begin{subequations}
\label{eq:ED}
\begin{align}
  \min \quad 
  & \sum_{t\in\mathcal{T}} \left( \sum_{i\in\mathcal{G}} C_i p_{it} + \sum_{j\in\mathcal{D}} C_j^d d_{jt} + \sum_{i\in\mathcal{I}} C_i^m m_{it} \right. \notag \\
  & \qquad \left. + \sum_{i\in\mathcal{W}} C_i^w w_{it} + \sum_{i\in\mathcal{R}} C_i^r r_{it} \right) \label{eq:ED:obj} \\
  \text{s.t.} \quad
  & \sum_{l\in\mathcal{L}_n^+} f_{lt} - \sum_{l\in\mathcal{L}_n^-} f_{lt} + \sum_{i\in\mathcal{G}_n} p_{it} + \sum_{i\in\mathcal{I}_n} (M_{it} - m_{it}) \notag \\
  & + \sum_{i\in\mathcal{W}_n} (W_{it} - w_{it}) + \sum_{i\in\mathcal{R}_n} (R_{it} - r_{it}) \notag \\
  & = \sum_{j\in\mathcal{D}_n} (D_{jt} - d_{jt}), \quad (\lambda_{nt}), \quad \forall n\in\mathcal{N}, t\in\mathcal{T}, \label{eq:ED:flow} \\
  & f_{lt} = B_l (\theta_{nt} - \theta_{mt}), \quad \forall l=(m,n)\in\mathcal{L}, t\in\mathcal{T}, \label{eq:ED:angle}\\ 
  & -RD_i \leq p_{it} - p_{i,t-1} \leq RU_i, \quad \forall i\in\mathcal{G}, t\in\mathcal{T}, \label{eq:ED:ramp} \\
  & -F_l^{max} \leq f_{lt} \leq F_l^{max}, \quad \forall l\in\mathcal{L}, t\in\mathcal{T}, \label{eq:ED:linecap}\\
  & \Theta_n^{min} \leq \theta_{nt} \leq \Theta_n^{max} \quad \forall n\in\mathcal{N}, t\in\mathcal{T}, \label{eq:ED:anglerange} \\
  & 0 \leq p_{it} \leq P_i^{max}, \quad \forall i\in\mathcal{G}, t\in\mathcal{T}, \label{eq:ED:maxgen}\\
  & 0 \leq d_{jt} \leq D_{jt}, \quad \forall j\in\mathcal{D}, t\in\mathcal{T}, \label{eq:ED:maxshed} \\
  & 0 \leq m_{it} \leq M_{jt}, \quad \forall i\in\mathcal{I}, t\in\mathcal{T}, \label{eq:ED:maximport} \\
  & 0 \leq w_{it} \leq W_{jt}, \quad \forall i\in\mathcal{W}, t\in\mathcal{T}, \label{eq:ED:maxwind} \\
  & 0 \leq r_{it} \leq R_{jt}, \quad \forall i\in\mathcal{R}, t\in\mathcal{T}. \label{eq:ED:maxrenewable}
\end{align}
\end{subequations}

Note that power is supplied from imports, (nonwind) renewables (e.g., bio-,
hydro- and geo-) and wind generations as well as from conventional thermal
generation units. Considering imports and nonwind renewables (we refer to these
simply as renewables in the following discussion) is important because they
account for a significant portion of the power generation in some systems. In
the CAISO system, for instance, imports and renewables account for 27\% and
25\% of the total system generation, respectively. Moreover, the analysis that
we present later indicates that dispatchable loads can reduce spillage of
imports and nonwind renewables. In the presented model we assume that imports as
well as renewable and wind power suppliers are not competitive agents in the
market (they are high-priority suppliers). Consequently, their supplies are
considered as negative demands for which we seek to minimize spillages at costs
$C_i^m,C_i^w$, and $C_i^r$. We also allow for load shedding at certain nodes at 
cost $C_j^d$,  which is set to the value of lost load (VOLL). 

The objective function \eqref{eq:ED:obj} is to minimize the {\em total dispatch
cost}: the sum of supply cost from conventional thermal generators, cost of the
load shedding, cost of the import spillage, cost of the wind power spillage,
and cost of the nonwind renewable spillage. Note that this objective can be
interpreted as maximizing social welfare, as defined in electricity market
clearing models (e.g., \cite{zavala2015stochastic}).  Equation
\eqref{eq:ED:flow} enforces the power flow balance of the network. Equation
\eqref{eq:ED:angle} represents a lossless model of power flow equations that
determines the power flow of line $l$ by the phase angle difference between two
buses $m$ and $n$. Equation \eqref{eq:ED:ramp} represents the ramping
constraints limiting the rate of change of generation levels. Constraints
\eqref{eq:ED:linecap} and \eqref{eq:ED:anglerange} represent the transmission
line capacity and the feasible phase angle range, respectively. Constraint
\eqref{eq:ED:maxgen} represents the generation capacity, and constraint
\eqref{eq:ED:maxshed} is a bound for the unserved load. Constraints
\eqref{eq:ED:maximport}-\eqref{eq:ED:maxrenewable} bound spillages of imports,
wind, and renewable supply, respectively.

\subsection{Optimal Placement of Dispatchable Loads}

We extend the ED model to account for optimal placement (OP) of dispatchable
loads at locations minimizing the expected total dispatch cost over multiple
wind and load scenarios. This OP model is cast as a two-stage stochastic
integer program. Following the convention in the literature (e.g.,
\cite{yu2014risk}) we provide the first-stage problem, the second-stage
problem, and the deterministic equivalent problem.

\subsubsection{First-Stage Problem}

The first-stage problem is given by
\begin{subequations}
\label{eq:OP1}
\begin{align}
  \min \; 
  & \mathbb{E} \left[ Q(x,\omega) \right] \label{eq:OP1:obj} \\
  \text{s.t.} \;
  & \sum_{n\in\mathcal{N}} x_n \leq K, \label{eq:OP1:const1}\\
  & x_n \geq 0, \text{ integer} \quad \forall n\in\mathcal{N},\label{eq:OP1:const2}
\end{align}
\end{subequations}
where $Q(x,\omega)$ is the recourse function of the first-stage variable $x$ for a given scenario $\omega$. The first-stage decision variable $x_n$ is a here-and-now decision to represent the number and locations of dispatchable loads to be installed at location $n\in\mathcal{N}$. Equation \eqref{eq:OP1:const1} is a budget constraint for dispatchable loads.

\subsubsection{Second-Stage Problem}

The second-stage problem is given as the recourse function $Q(x,\omega)$
defined for each scenario $\omega\in\Omega$ as follows:
\begin{subequations}
\label{eq:OP2}
\begin{align}
  \min \; 
  & \sum_{t\in\mathcal{T}} \left( \sum_{i\in\mathcal{G}} C_i p_{it}(\omega) + \sum_{j\in\mathcal{D}} C_j^d d_{jt}(\omega) + \sum_{i\in\mathcal{I}} C_i^m m_{it}(\omega) \right. \notag \\
  & \qquad \left. + \sum_{i\in\mathcal{W}} C_i^w w_{it}(\omega) + \sum_{i\in\mathcal{R}} C_i^r r_{it}(\omega) \right)  \label{eq:OP2:obj} \\
  \text{s.t.} \;
  & 0 \leq u_{nt}(\omega) \leq U x_n, \quad \forall n\in\mathcal{N}, t\in\mathcal{T}, \label{eq:OP2:const3} \\
  & \sum_{l\in\mathcal{L}_n^+} f_{lt}(\omega) - \sum_{l\in\mathcal{L}_n^-} f_{lt}(\omega) + \sum_{i\in\mathcal{G}_n} p_{it}(\omega) \notag\\
  & + \sum_{i\in\mathcal{I}_n} (M_{it}-m_{it}(\omega)) + \sum_{i\in\mathcal{W}_n} (W_{it}(\omega)-w_{it}(\omega)) \notag\\
  & + \sum_{i\in\mathcal{R}_n} (R_{it}-r_{it}(\omega)) - u_{nt}(\omega) \notag\\
  & = \sum_{j\in\mathcal{D}_n} (D_{jt} - d_{jt}(\omega)), \quad \forall n\in\mathcal{N}, t\in\mathcal{T}, \label{eq:OP2:const1} \\
  & f_{lt}(\omega) = B_l (\theta_{nt}(\omega) - \theta_{mt}(\omega)), \notag\\
  & \qquad \forall l=(m,n)\in\mathcal{L}, t\in\mathcal{T}, \label{eq:OP2:angle}\\ 
  & -RD_i \leq p_{it}(\omega) - p_{i,t-1,}(\omega) \leq RU_i, \quad \forall i\in\mathcal{G}, t\in\mathcal{T}, \label{eq:OP2:ramp} \\
  & -F_l^{max} \leq f_{lt}(\omega) \leq F_l^{max}, \quad \forall l\in\mathcal{L}, t\in\mathcal{T}, \label{eq:OP2:linecap}\\
  & \Theta_{n}^{min} \leq \theta_{nt}(\omega) \leq \Theta_n^{max} \quad \forall n\in\mathcal{N}, t\in\mathcal{T}, \label{eq:OP2:anglerange} \\
  & 0 \leq p_{it}(\omega) \leq P_i^{max}, \quad \forall i\in\mathcal{G}, t\in\mathcal{T}, \label{eq:OP2:maxgen}\\
  & 0 \leq d_{jt}(\omega) \leq D_{jt}, \quad \forall j\in\mathcal{D}, t\in\mathcal{T}, \label{eq:OP2:maxshed} \\
  & 0 \leq m_{it}(\omega) \leq M_{jt}, \quad \forall i\in\mathcal{I}, t\in\mathcal{T}, \label{eq:OP2:maximport} \\
  & 0 \leq w_{it}(\omega) \leq W_{jt}, \quad \forall i\in\mathcal{W}, t\in\mathcal{T}, \label{eq:OP2:maxwind} \\
  & 0 \leq r_{it}(\omega) \leq R_{jt}, \quad \forall i\in\mathcal{R}, t\in\mathcal{T}. \label{eq:OP2:maxrenewable}
\end{align}
\end{subequations}

The objective function \eqref{eq:OP2:obj} is the total dispatch cost with
uncertainty arising from wind supply scenarios. The second-stage decisions
involve flows, angles, supply, loads, and spillages. Equation
\eqref{eq:OP2:const3} is a capacity constraint for dispatchable loads. The
capacity of a dispatchable load is given by $U$. We note that constraints
\eqref{eq:OP2:const1} include dispatchable loads.

\subsubsection{Deterministic Equivalent Problem}

Assuming that $\omega$ has finite support on $\Omega$, the two-stage stochastic
programming problem in \eqref{eq:OP1} and \eqref{eq:OP2} can be formulated as a
deterministic equivalent problem:
\begin{subequations}
\label{eq:OP3}
\begin{align}
  \min \; 
  & \sum_{\omega\in\Omega} \pi(\omega) \sum_{t\in\mathcal{T}} \left( \sum_{i\in\mathcal{G}} C_i p_{it}(\omega) + \sum_{j\in\mathcal{D}} C_j^d d_{jt}(\omega) \right. \notag \\
  & \qquad + \sum_{i\in\mathcal{I}} C_i^m m_{it}(\omega) + \sum_{i\in\mathcal{W}} C_i^w w_{it}(\omega) \notag \\
  & \qquad \left. + \sum_{i\in\mathcal{R}} C_i^r r_{it}(\omega) \right)  \label{eq:OP3:obj} \\
  \text{s.t.} \;
  & \eqref{eq:OP1:const1}, \eqref{eq:OP1:const2}, \notag \\
  & \eqref{eq:OP2:const3}-\eqref{eq:OP2:maxrenewable}, \quad \forall \omega\in\Omega. \notag
\end{align}
\end{subequations}

\subsection{Performance Metrics}

We define metrics of interest for our analysis. The total amount of power
dispatched (absorbed into the system) is given by
\begin{align}
  P_{ED}^{Dispatch} := & \sum_{t\in\mathcal{T}} \left[ \sum_{i\in\mathcal{G}} p_{it} + \sum_{i\in\mathcal{I}} ( M_{it} - m_{it} ) \right. \notag \\
  & \quad \left. + \sum_{i\in\mathcal{W}} ( W_{it} - w_{it} ) + \sum_{i\in\mathcal{R}} ( R_{it} - r_{it} ) \right].
\end{align}
We differentiate the total power absorbed at positive LMPs and nonpositive
LMPs. To do so, we define $\mathcal{N}_t^+ := \{n\in\mathcal{N} : \lambda_{nt}
> 0\}$, where $\lambda_{nt}$ is an optimal dual variable value of the ED model.
The total power absorbed at positive LMPs and nonpositive LMPs for ED is given
respectively by
\begin{align}
  P_{ED}^{LMP>0} := 
  &\sum_{t\in\mathcal{T}} \sum_{n \in \mathcal{N}_t^+} \left[ \sum_{i\in\mathcal{G}_n} p_{it} + \sum_{i\in\mathcal{I}_n} \left( M_{it} - m_{it} \right) \right. \notag\\
  &\left .+ \sum_{i\in\mathcal{W}_n} \left(W_{it} - w_{it} \right) + \sum_{i\in\mathcal{R}_n} \left(R_{it}-r_{it}\right)\right]\label{eq:ED:poslmp}
\end{align}
and $P_{ED}^{LMP\leq 0} := P_{ED}^{Dispatch} - P_{ED}^{LMP>0}$. 

We define stranded power as $P_{ED}^{LMP\leq 0} + P_{ED}^{Spillage}$. Wind
penetration levels (\%) are defined as $100\times \sum_{i\in\mathcal{W}}
\sum_{t\in\mathcal{T}} W_{it} / \sum_{j\in\mathcal{D}} \sum_{t\in\mathcal{T}}
D_{jt}$.
The RPS is defined by the ratio of the renewable power absorbed to the total
loads:
\begin{align}
\label{eq:RPS}
 100\times  \frac{\sum_{t\in\mathcal{T}} \left[\sum_{i\in\mathcal{R}} (R_{it}-r_{it}) + \sum_{i\in\mathcal{W}} (W_{it} - w_{it}) \right]}{\sum_{j\in\mathcal{D}} \sum_{t\in\mathcal{T}} D_{jt}}.
\end{align}
For the OP model the metrics are defined for each $s\in \mathcal{S}$ from which
we compute expected values. 

The achieved capacity (\%) at time $t$ for a given scenario $s$ is defined as
the ratio of the number of data centers served by a positive stranded power:
$\sum_{n\in\mathcal{N}} \mathbf{1}(u_{nts}) / K$, where $\mathbf{1}(u_{nts}) =
1$ if $u_{nts} > 0$ and zero otherwise.


\section{Computational Results and Analyses}
\label{section:experiments}

\begin{figure}
\centering
\subfloat[Net load for each day type.]{
  \includegraphics[width=0.23\textwidth]{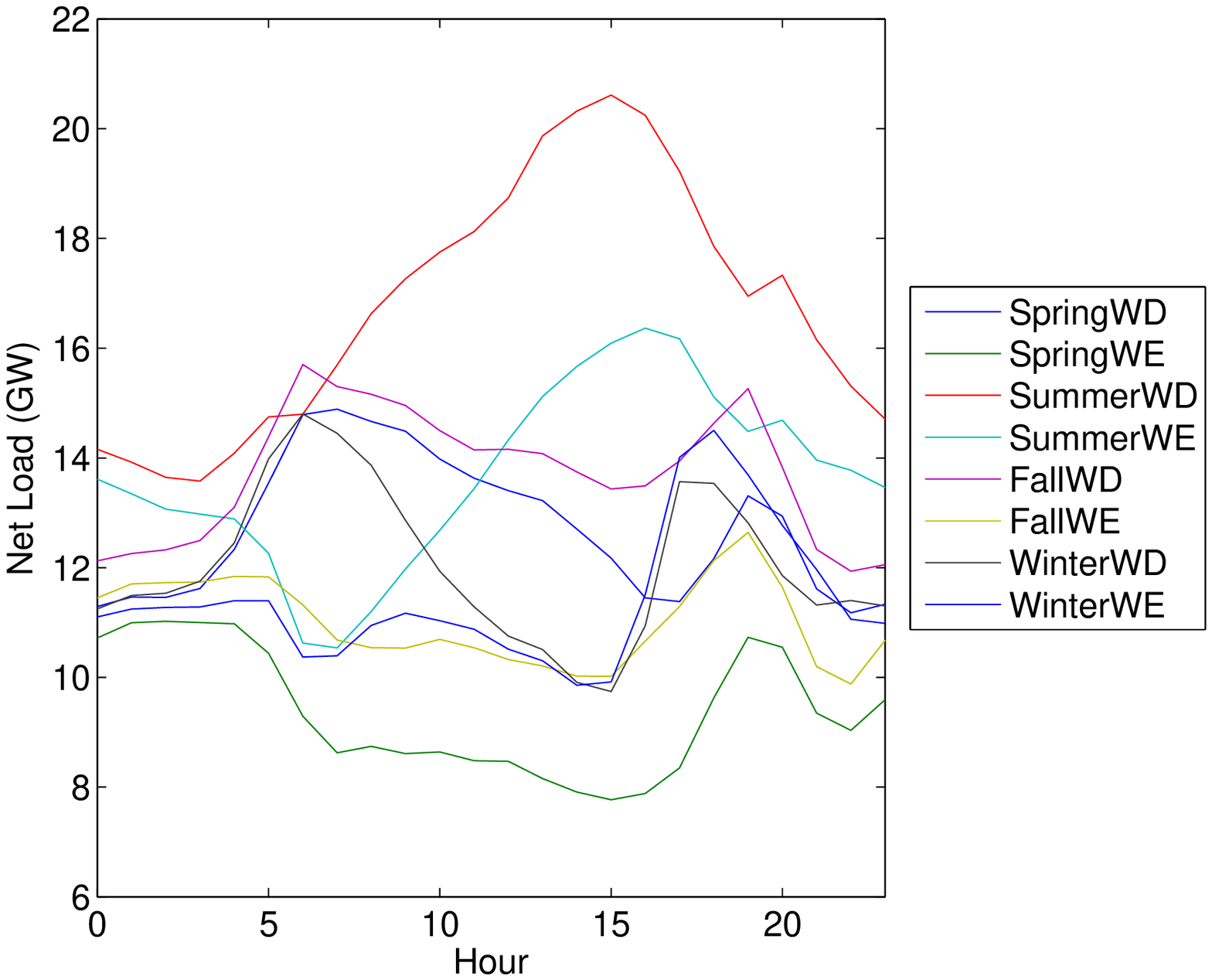}\label{fig:netload}
}
\centering
\subfloat[Wind power supply (GW) at 15\% level during summer day.]{
  \includegraphics[width=0.23\textwidth]{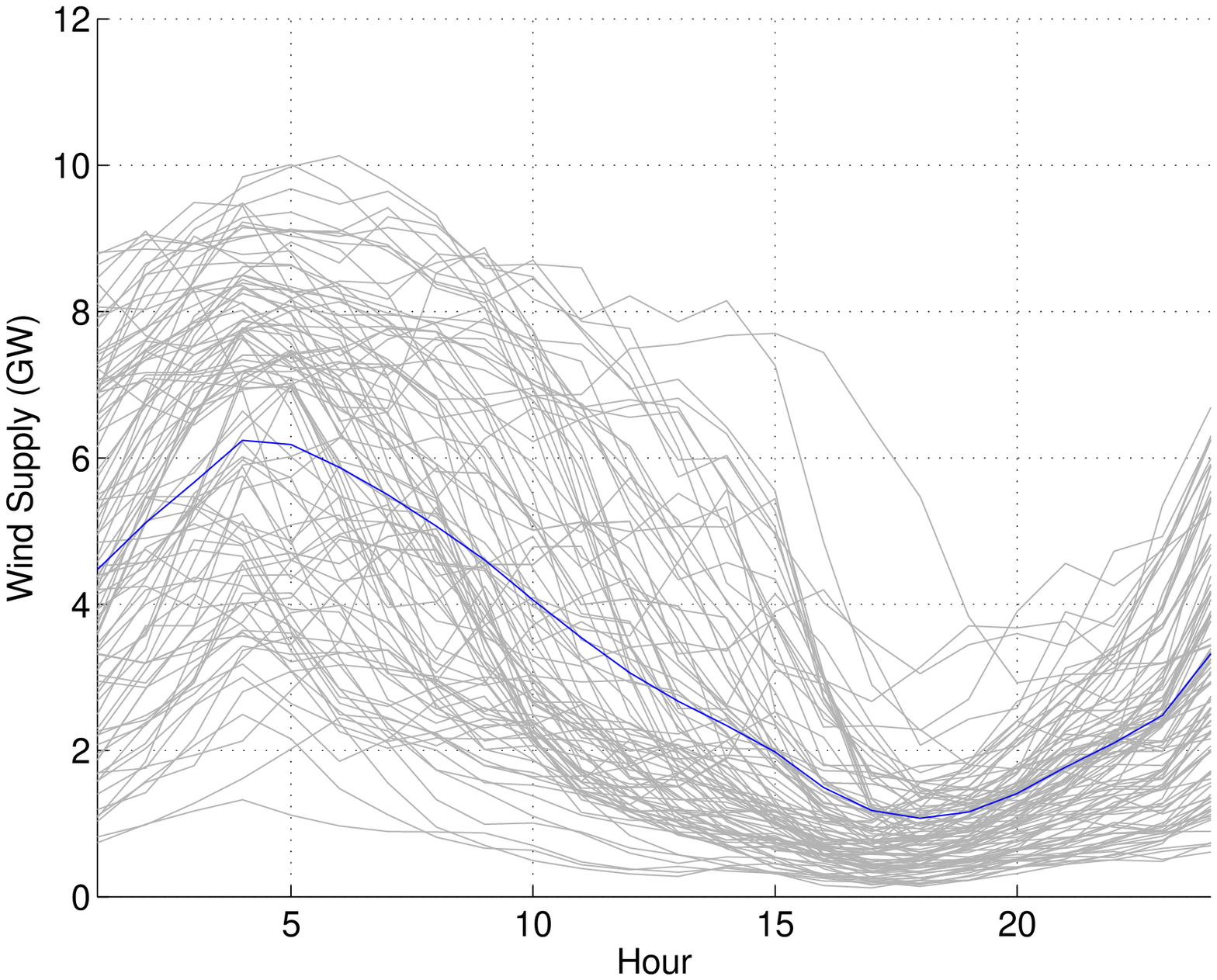}\label{fig:windprod}
}
\caption{Net load profile and wind power supply.}
\end{figure}

We study a test system of CAISO interconnected with the WECC. The system
consists of 225 buses, 375 transmission lines, 130 generation units, 40 loads,
and 5 wind power generation units. We consider a 24-hour horizon with hourly
intervals. We use network topology, import supply, renewables, wind production,
and load data from \cite{papavasiliou2013multiarea}. Imports flow into the
system through 5 boundary buses. Renewable power is generated from biogas,
hydrothermal, and geothermal generators at 11 buses.  For this system, the
generation capacity that excludes imports, renewables, and wind power is 31.2
GW. We consider load profiles for 8 day types: spring, summer, fall, and
winter; and weekday (WD) and weekend (WE). Table~\ref{tab:avgload}
reports average loads, imports, renewables, and net loads (load minus imports,
renewables, and wind supply).  Figure~\ref{fig:netload} shows the net loads for
the different day types.

\begin{table}
\centering
\caption{Average load, imports, renewable supply, and net load (MW)}
\label{tab:avgload}
\vspace*{-0.1in}
\begin{tabular}{lcccc}
  \hline
  & Load & Imports & Renewable & Net Load\\
  \hline
  SpringWD & 26,868 & 7,478 & 6,681 & 12,708 \\
  SpringWE & 23,980 & 7,608 & 6,998 &  9,373 \\
  SummerWD & 31,089 & 7,678 & 6,672 & 16,737 \\
  SummerWE & 28,184 & 7,400 & 7,124 & 13,659 \\
  FallWD   & 28,055 & 7,675 & 6,657 & 13,722 \\
  FallWE   & 25,186 & 7,108 & 7,065 & 11,012 \\
  WinterWD & 26,352 & 7,663 & 6,634 & 12,054 \\
  WinterWE & 23,708 & 6,800 & 5,581 & 11,399 \\
  \hline
\end{tabular}
\end{table}

For each day type, we use the 1,000 specific wind power production scenarios
taken from \cite{papavasiliou2013multiarea} where wind power scenarios are 15\%
of load, representing the 2020 RPS target of California
\cite{papavasiliou2013multiarea} as illustrated in Figure~\ref{fig:windprod}.
Each scenario has the same probability $\pi(\omega)=0.001$ for
$\omega\in\Omega$.  We chose 1,000 scenarios because they are enough for
providing the mean objective values with statistical significance p-value $<
0.0001$ for the four cases. However, the number of sample sizes is not known a
priori in general.  For readability, we plot only 100 scenarios for wind power
(grey lines) to highlight the variability, and we also plot the corresponding
mean (blue line).  The same wind scenarios are used for weekdays and weekends
of the same season.  We explore a range of additional wind penetration levels:
5\%, 15\%, 30\%, and 50\% of load,

As is typical due to longer-term commitments and with the goal of reducing carbon
emissions, we assume that imports and renewables are higher priority (i.e.,
nondispatchable) but can be spilled, if necessary, at a cost of 1,000 and 2,000
\$/MWh, respectively. We use a load-shedding cost (VOLL)  of 1,000 \$/MWh and a
wind spillage cost of 100 \$/MWh. These values are typical for ISO settings
and are chosen to impose a relative priority on different products.  The twenty
additional 200 MW data centers reflect the rapid growth expected in cloud
computing over the next few decades, particularly in the western region.  ICT
power consumption is estimated at 8\% today and projected to grow to as much
as 4$\sim$10\% per decade \cite{GreenPeaceCC15}.

We analyze four cases.
\begin{itemize}
  \item Case 1: Base, WECC configuration as described above.

  \item Case 2: Case 1 with twenty additional 200 MW data centers that total 4
GW additional load (96 GWh per day).  Each data center is a continuous 200 MW
load and subject to VOLL penalties.  Data-center locations were chosen
arbitrarily to reflect choices driven by external business considerations
(e.g., networking, proximity to customers, and geographic diversity). The data
center loads are inflexible (non-dispatchable).

  \item Case 3: Case 2 with collocated wind farms at each data center (with
non-dispatchable loads), sized to match total load over 12 months.  Because of the
typical wind capacity factor of 30\%, the peak generation of these farms is
typically 3x greater than the 200 MW data-center load.

  \item Case 4: Case 1 with twenty additional 200 MW data centers operated as
dispatchable loads.  The ISO determines the power consumption of each
dispatchable load each hour at no penalty cost.  The data centers are
positioned to minimize overall system dispatch cost across all wind and load
scenarios by solving the OP model \eqref{eq:OP1}.

\end{itemize}

In Cases 1, 2, and 3, we solve the ED model \eqref{eq:ED}, minimizing the total
dispatch cost for all wind and load scenarios as we increase wind levels. In
Case 4, the ED is subsumed within the OP model, and the solution minimizes the
same metrics in addition to optimal placement of dispatchable loads.

The cases vary in total generation and load.  We compare percentages relative
to Case 1 (i.e., the base system).  For example, the loads in Case 2 and 3 are
higher because of the addition of data centers with non-dispatchable loads, and
Case 4 falls in the middle because of its variable dispatch capacity.  We use
consistent wind penetration numbers, ignoring the additional wind generation in
Case 3.  The simple treatment of loads affects ``real'' wind penetration only
1$\sim$6\%, far smaller than the resulting spillage.  Excluding data-center
wind power in Case 3 produces conservative estimates of spillage and stranded
power, painting Case 3 in the most favorable light.

\begin{figure}
  \centering
  \includegraphics[width=0.25\textwidth]{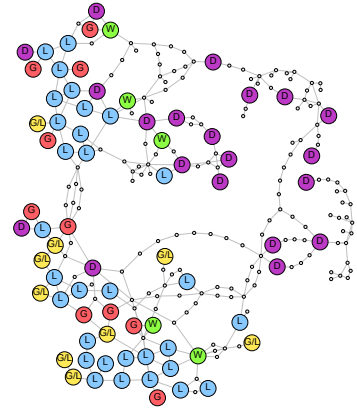}
  \caption{Node-edge network representation of test system (Cases 2 and 3). The 
buses with generation units, loads, both generation units and loads, wind units, 
and data centers are labeled as G, L, G/L, W and D, respectively.}
  \label{fig:wecc}
\end{figure}

Figure~\ref{fig:wecc} presents a node-edge network representation of the test
system with the twenty data center locations of Case 2. We note that this
network does not represent the actual geographical locations of the actual
buses and the lines of the system. The network was generated by using the {\tt
gephi} package \cite{bastian2009gephi}. The network system under study is
large, and we evaluate a large number of scenarios and system configurations.
Thus, the analysis performed is computationally intensive. Cases 1, 2, and 3
solve the ED model \eqref{eq:ED} for 1,000 wind-power scenarios and for each
day type and season. These represent a total of 32,008 linear programs (LPs)
for each case. Each LP has 19,008 continuous variables and 17,544 linear
constraints. Case 4 solves the OP model \eqref{eq:OP1} including all the 1,000
wind scenarios; each OP instance is solved for the different day type and
season and different wind power levels. We thus solve a total of 24 large-scale
stochastic mixed-integer programs (for nonzero wind levels) and 8 deterministic
mixed-integer programs (at zero wind level). Each {\em OP model has 225 general
integer variables, 24,408,000 continuous variables, and 22,944,001
constraints}. The ED model was implemented in JuMP \cite{lubin2015computing}
and the OP model in StochJuMP \cite{huchette2014parallel}. ED is solved with
CPLEX 12.6.1, and OP is solved by using the parallel Benders decomposition
implementation of the open-source package DSP \cite{kim2015algorithmic}.  DSP
is a high-performance optimization package capable of targeting large-scale
stochastic programming problems. 

All computations were performed on {\tt Blues}, a 310-node computing cluster at
Argonne National Laboratory. Each computing node has two octo-core 2.6 GHz Xeon
processors and 64 GB of RAM. The computational experiments required the
solution of 128,032 linear programs, 8 mixed-integer programs, and 24 stochastic
mixed-integer programs and post-processing of solution data. In addition to
the parallel solver DSP, we used SWIFT \cite{wilde2011SWIFT}, a script language
for dynamically allocating computing jobs in the parallel cluster. In
particular, SWIFT allows us to solve multiple optimization problems in
parallel, while each problem is in turn solved in parallel by using DSP. This
setting enables us to perform computationally intensive experiments. To give an
idea of the efficiency achieved, a single stochastic optimization problem OP
was solved in 10 minutes using 1,000 parallel cores (the same problem would
have required 166 hours if run in a standard serial computer).  Running the
entire set of computational experiments using DSP and SWIFT required only 25
hours of wall-clock time with 2,000 computing cores (50,000 core-hours).  In
contrast, performing such experiments in serial would have required nearly 6
months. 

The results reveal important trends that we summarize below. We then present
numerical results to illustrate these trends.

\begin{itemize}

\item Case 1 reveals significant spillage and stranded power in the base system
due to imports and nonwind renewables.  This finding is consistent with the
observations of Section \ref{sec:stranded}. We also observe that, as expected,
dispatch cost decreases initially as cheaper wind power displaces thermal
generation, but eventually increases because of stranded power penalties. We
also see that the variance of the cost increases dramatically as wind level is
increased, indicating that the system becomes more vulnerable to wind power
variations.

\item Case 2 reveals that positioning large data centers decreases system cost,
even if locations are chosen arbitrarily and loads are inflexible. The reason
is that the loads put stranded power to work, reducing penalties and moderately
reducing system cost.  We also observe that while cost is decreased, the
variance of the cost is not improved (compared with Case 1).

\item Case 3 reveals that collocating data centers at wind-farm locations gives
little benefit to system cost. The slight benefit comes from wind power used to
offset the data-center loads, but stranded power is increased as a result of
load inflexibility.  Case 3 also reveals that collocation of data centers and
wind farms does not reduce carbon emissions. In fact, it increases thermal
generation (and thus emissions) because of increasing stranded power.  We found
that system cost variance was reduced significantly compared with that of Cases
1 and 2. We attribute this reduction to better utilization of stranded power at
data-center locations.  This result thus highlights that stranded power can
affect system vulnerability.

\item Case 4 reveals that optimally positioned dispatchable loads can reduce
power spillage from all sources (imports, nonwind renewables, and wind), not
just wind (as in Case 3). Strategic positioning results in decreased system
cost and decreased use of thermal generation (and thus emissions). We also find
that cost variance is dramatically reduced compared with all other cases,
indicating that the system can better manage wind power fluctuations due to
increased network flow control.  Case 4 also reveals that up to 60$\sim$80\% of
data center capacity can be achieved at high wind-penetration levels.  This
result occurs even if loads can be adjusted at no cost. Consequently,
data-center owners have a natural economic incentive to provide flexibility.
The incentives are provided by better stranded power utilization.

\end{itemize}

\begin{figure}[!htb]
  \centering
  \includegraphics[width=0.45\textwidth]{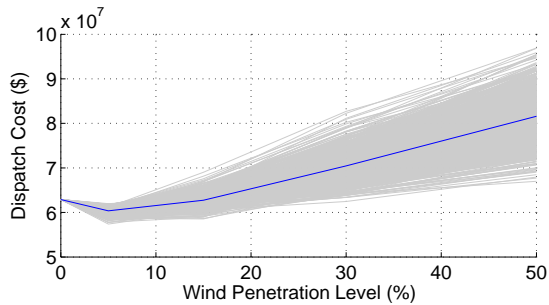}
  \caption{Dispatch cost at different wind levels (Case 1).}
  \label{fig:social_welfare}
\end{figure}

\subsection{Base System}

We first analyze the impact of increasing wind levels in Case 1.
Figure~\ref{fig:social_welfare} shows that the average daily dispatch cost
decreases below a 5\% wind level, because of the use of cheap wind power. As
wind levels increase, however, the dispatch cost is increased by 26\% ($\$2.1$
MUSD/day) relative to the system at 0\% wind level. This increase is a combined
effect of using more thermal generation to account for wind power variability
and penalties induced by power spillage. Moreover, the system cost becomes more
variable as we increase wind levels. In particular, the standard deviation is
$\$5.2$ MUSD/day at a 50\% wind level and $\$0.7$ MUSD/day at a 5\% wind level.
This variability indicates that the system becomes more vulnerable at high wind
levels.

\begin{figure}
  \centering
  \subfloat[Absorbed Power]{
    \includegraphics[width=0.23\textwidth]{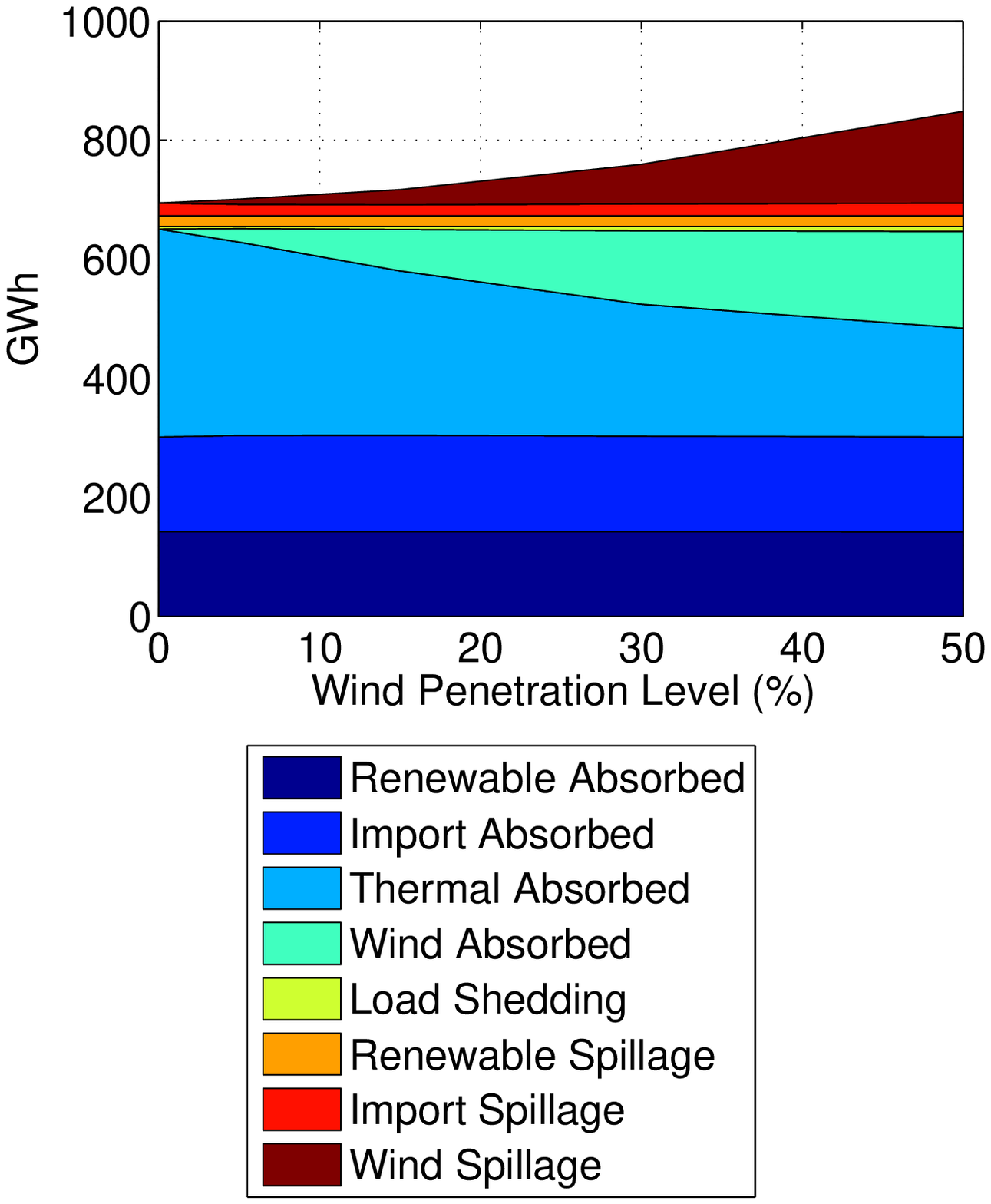}
  }
  \subfloat[Total Supply]{
    \includegraphics[width=0.23\textwidth]{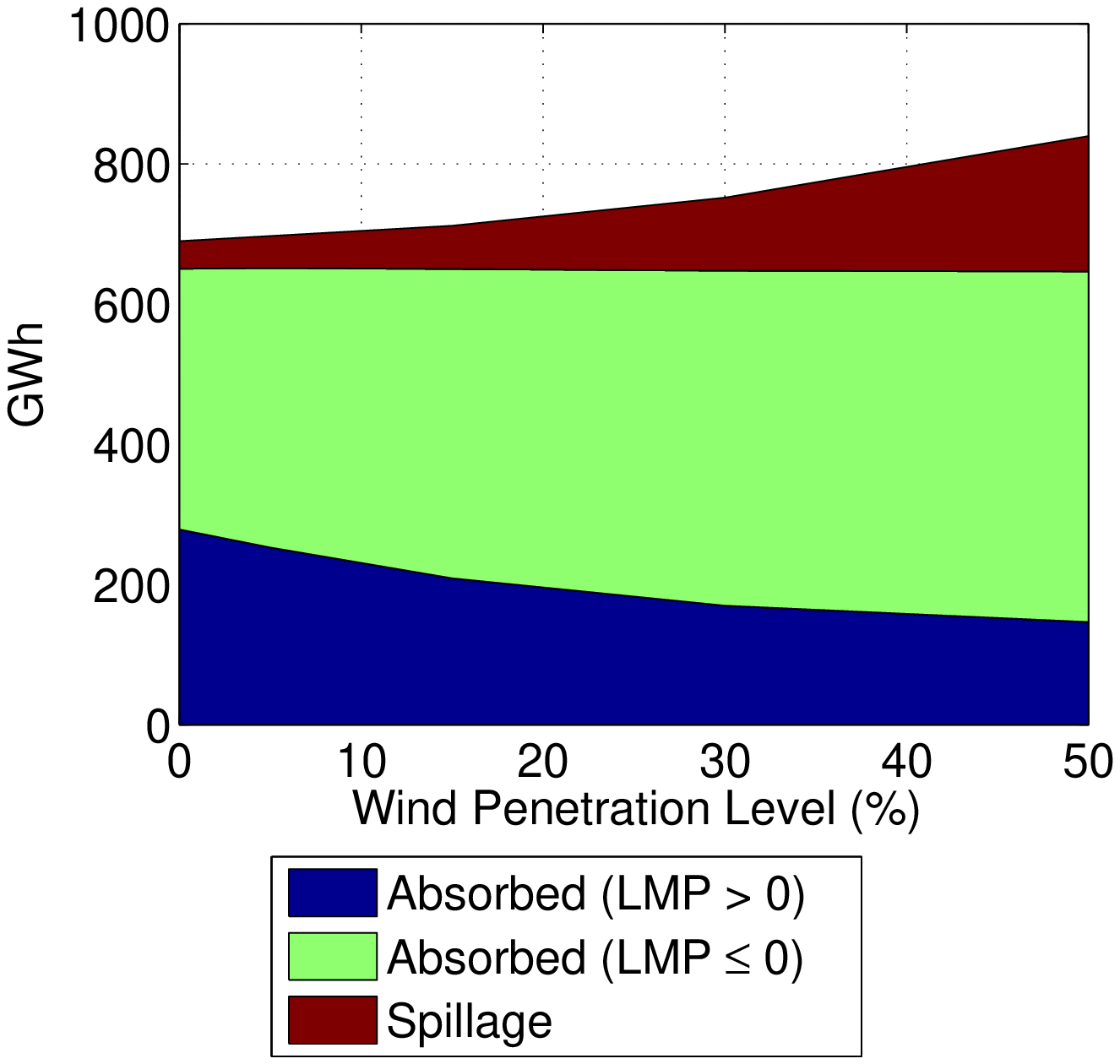}
  }
  \caption{Absorbed power, spillage, and total supply for increasing wind level (Case 1).}
  \label{fig:absorbed-base}
\end{figure}

Figure~\ref{fig:absorbed-base}a shows daily average power and spillage by
generation source.  As wind penetration is increased form 0\% to 50\%, thermal
generation decreases by 48\% (from 349 to 182 GWh).  To absorb the variability
of the wind, total supply increases by 21\% (from 689 to 839 GWh), even if
there is increase in the load.  At 50\% wind penetration level, 23\% of
generation is spilled ({\em not absorbed} into the system).  This result
reflects the difficulty in achieving high RPS because of wind variability.
Daily load shed for this case is also 8 GWh. 

Figure~\ref{fig:absorbed-base}b shows daily average power by LMP value (price)
and spillages.  We recall that stranded power is the sum of spillage and power
absorbed into the system at negative price (LMP$\leq0$). We observe that as
wind levels increase, both spillage and stranded power increase.  Stranded
power increases from 60\% at a 0\% wind level to 83\% at a 50\% wind level.
While the total economic return for a generator may not be reflected in LMP
alone, we note that the amount of power absorbed at positive price (profitable
power) does not increase after a 15\% wind level (see
Figure~\ref{fig:WindByEconomicOthers}b).

\subsection{Adding Data Centers}

\begin{table}
\centering
\caption{Thermal supply and spillages at different wind levels (GWh)}
\label{tab:dispatch}
\begin{tabular}{lrrrrrr}
  \hline
        & Wind  &     & Thermal & Wind      & Import   & Renewable \\
        & Level & RPS & Supply  & Spillage  & Spillage & Spillage \\
  \hline
  Case 1 & 0\%  & 22\% & 349 & 0   & 21 & 18 \\
         & 5\%  & 25\% & 325 & 8   & 19 & 18 \\
         & 15\% & 32\% & 276 & 25  & 19 & 18 \\
         & 30\% & 41\% & 221 & 66  & 20 & 18 \\
         & 50\% & 47\% & 182 & 154 & 21 & 18 \\
  \hline
  Case 2 & 0\%  & 22\% & 439 & 0   & 14 & 18 \\
         & 5\%  & 25\% & 413 & 6   & 12 & 18 \\
         & 15\% & 33\% & 361 & 19  & 12 & 18 \\
         & 30\% & 43\% & 292 & 48  & 12 & 18 \\
         & 50\% & 51\% & 238 & 122 & 13 & 18 \\
  \hline
  Case 3 & 0\%  & 22\% & 374 & 8   & 11 & 17 \\ 
         & 5\%  & 25\% & 355 & 17  & 11 & 17 \\ 
         & 15\% & 32\% & 315 & 39  & 11 & 17 \\ 
         & 30\% & 40\% & 272 & 86  & 11 & 17 \\ 
         & 50\% & 45\% & 237 & 176 & 11 & 17 \\ 
  \hline
  Case 4 & 0\%  & 24\% & 358 & 0   & 0 & 0 \\
         & 5\%  & 28\% & 362 & 2   & 9 & 5 \\
         & 15\% & 36\% & 309 & 14  & 11 & 7 \\
         & 30\% & 45\% & 255 & 49  & 11 & 7 \\
         & 50\% & 52\% & 220 & 130 & 11 & 8 \\
  \hline
\end{tabular}
\end{table}

\begin{figure}
  \centering
  \includegraphics[width=0.45\textwidth]{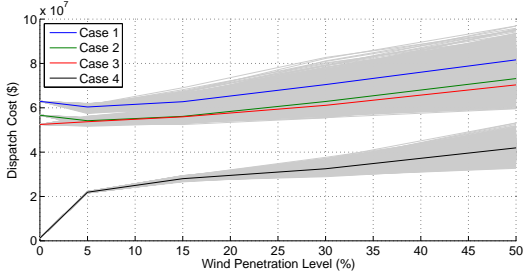}
  \caption{Dispatch cost at different wind levels.}
  \label{fig:social_welfare_comparison}
\end{figure}

Figure~\ref{fig:social_welfare_comparison} compares dispatch costs for all
cases. We note that the dispatch cost of the base system is decreased in all
cases. The dispatch costs of the base system are decreased in Cases 2 and 3 by
less than 5\% at 0\% wind level, respectively, and by less than 13\% at a 50\%
wind level. These results indicate that adding data centers has a beneficial
effect and that this value increases as more stranded power is present. Case 4
reduces the dispatch cost dramatically. A relative reduction of 98\% is
observed at 0\% wind level and a relative reduction of 49\% is observed at 50\%
wind level (compared with Case 1).  As seen in Table~\ref{tab:dispatch}, the
dramatic decrease in cost is due to minimization of spillage (penalized at
large values).  At a 0\% wind level, in particular, dispatchable loads fully
eliminate spillages. In Cases 2 and 3 we can see reductions in spillages, but
these are much smaller than those observed in Case 4. In particular, we note
that optimally placed dispatchable loads (Case 4) favor spillage reductions of
nonwind renewable supply over wind supply. This result supports the conclusion
that dispatchable loads, when strategically placed, can provide much greater
benefits for the power grid.  The reason is that optimal placement allows them
to eliminate spillages from a variety of generators in types and locations and
provides flexibility to manipulate network flows. 

In Figure~\ref{fig:social_welfare_comparison} we observe that significant value
is obtained in Case 4 at all wind levels, with benefits as great as 40\% at the
50\% wind level.  At low wind penetration, the dispatchable loads eliminate
essentially all spillage, dramatically reducing associated penalties.  As the
wind penetration level increases, the gap with respect to the base system is
reduced. The decreased benefit is due to the large amounts of spillage that are
introduced at high wind levels and that cannot be fully eliminated even with
optimally located dispatchable loads. This situation is observed in Table
\ref{tab:dispatch} where wind power spillage increases proportionally to the
wind level. Fully eliminating spillage would require additional data centers.

\begin{figure}
  \centering
  \includegraphics[width=0.45\textwidth]{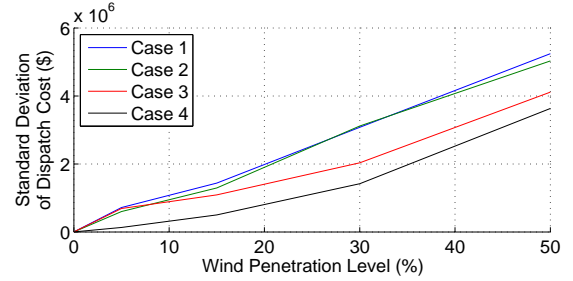}
  \caption{Standard deviation of dispatch costs at different wind levels.}
  \label{fig:social_welfare_variance}
\end{figure}

A surprising result is that the system cost variance is dramatically reduced
with dispatchable loads.  Figure \ref{fig:social_welfare_variance} illustrates
this. In particular, the standard deviation for Case 4 is $\$0.1$ MUSD/day at a
5\% wind level and $\$3.6$ MUSD/day at a 50\% wind level. We recall that the
standard deviations for Case 1 are $\$0.7$ MUSD/day at a 5\% wind level and
$\$5.2$ MUSD/day at a 50\% wind level. For Case 2 the standard deviation is
$\$0.6$ and $\$5$ MUSD/day at 5\% wind level and 50\% wind level, respectively.
For Case 3 the standard deviation is $\$0.6$ and $\$4.1$ MUSD/day at 5\% wind
level and 50\% wind level, respectively. For Case 4 we also note that variances
are negligible for wind levels below 10\% and remain small for wind levels
below 20\%. In contrast, the variances for Cases 1, 2, and 3 quickly increase
with the wind level. The reduction in cost variance is the result of additional
system spatiotemporal flexibility provided by dispatchable loads at different
locations, which allows for network flow control.

\subsection{Impact on Wind Generators}

\begin{figure}
  \centering
  \subfloat[Case 1]{
    \includegraphics[width=0.23\textwidth]{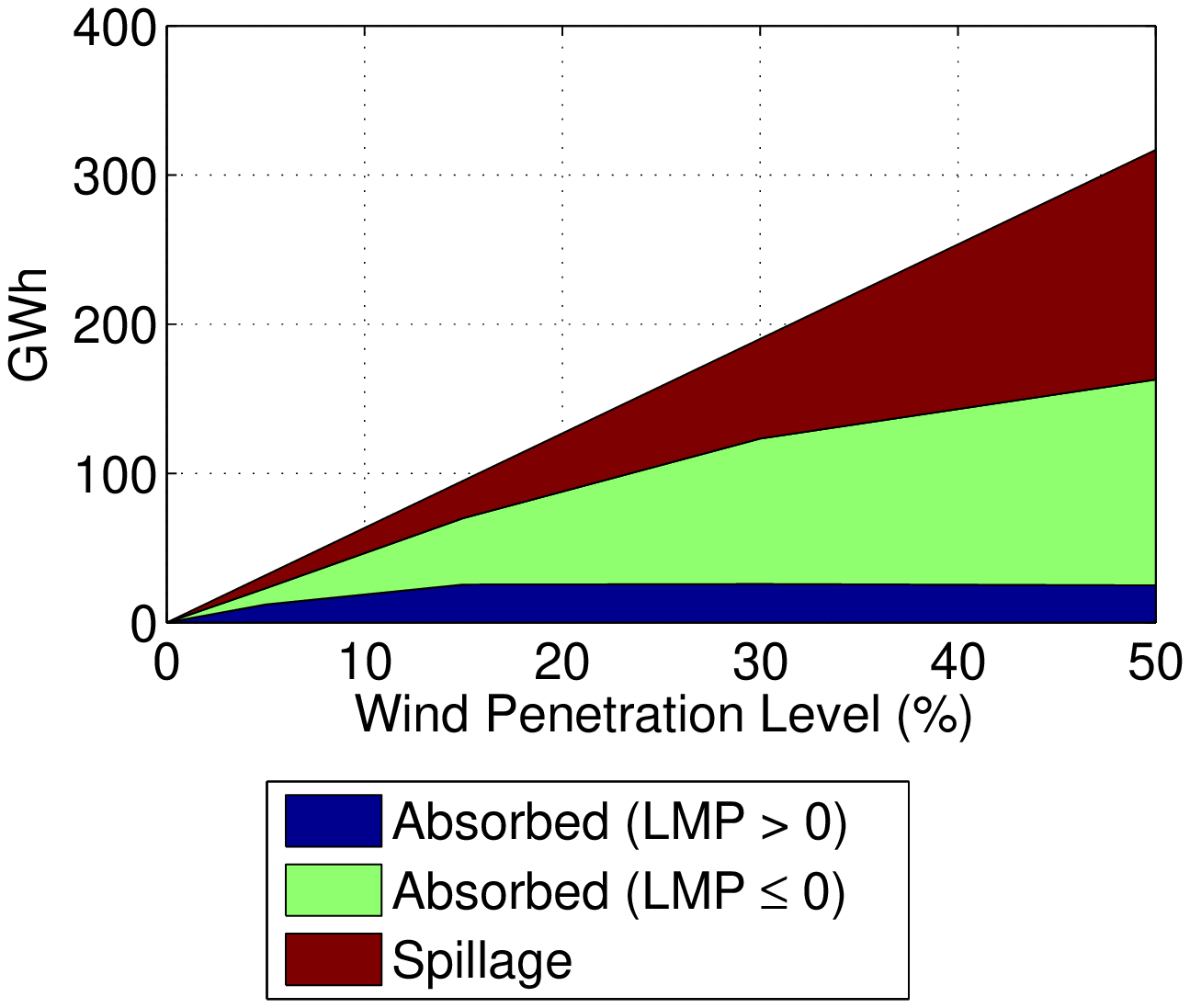}
  }
  \subfloat[Case 4]{
    \includegraphics[width=0.23\textwidth]{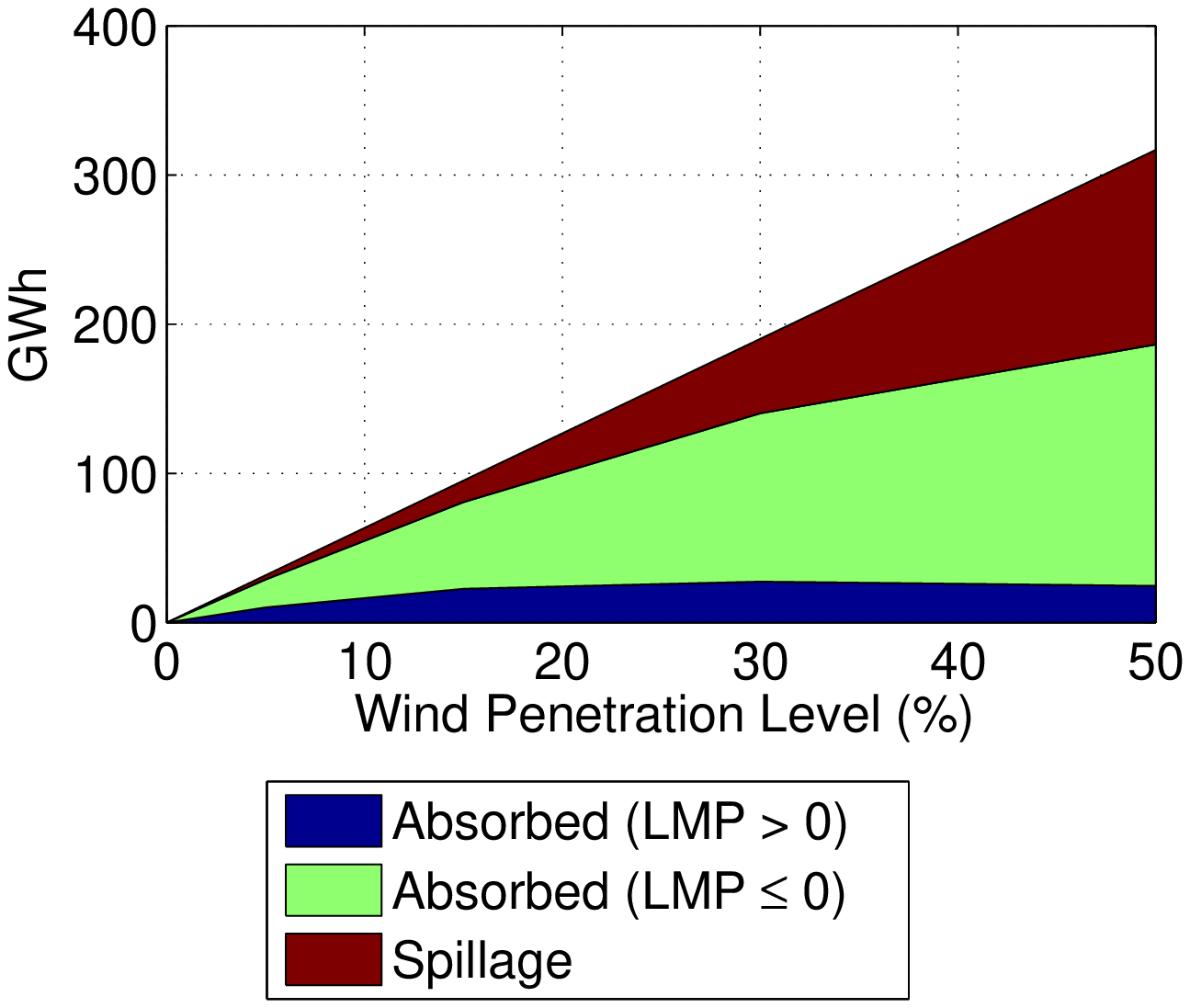}
  }
  \caption{Wind supply by LMP and spillage.}
  \label{fig:WindByEconomicOthers}
\end{figure}

Of particular interest are the generation, spillage, and uneconomic and
economic generation of wind power as the wind penetration level increases.
Figure \ref{fig:WindByEconomicOthers} shows the growing wind supply, what
fraction is spilled, and what fraction is absorbed in both uneconomic and
economic terms.  For brevity we show results only for Cases 1 and 4.  Adding
dispatchable loads decreases wind spillage significantly (to zero at low wind
penetration and by more than 15\% at high penetration).  However, there is a
complementary increase in the uneconomic power accepted by the grid, so the
total stranded power remains large.

\subsection{Data Center Achieved Capacity} 

\begin{figure}
  \centering
  \subfloat[Average achieved capacity daily time profile at different wind penetration levels.]{
    \includegraphics[width=0.45\textwidth]{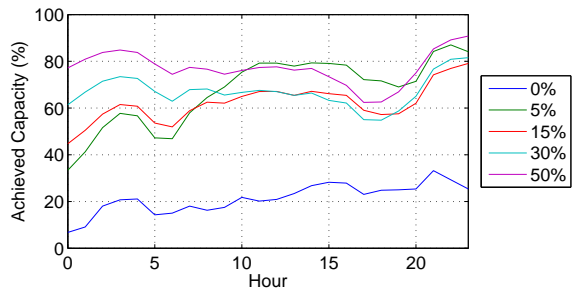}\label{fig:duty_factor1}
}\\
  \subfloat[Total achieved capacity increase in wind level]{
    \includegraphics[width=0.45\textwidth]{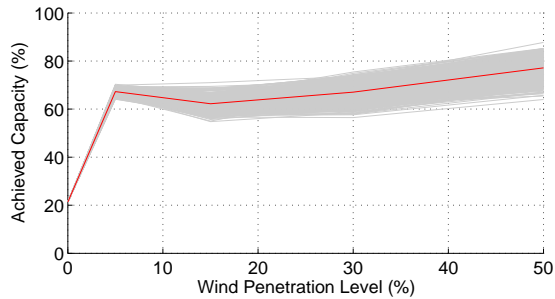}\label{fig:duty_factor2}
  }
  \caption{Achieved fraction of data-center capacity in Case 4.}
  \label{fig:duty_factor}
\end{figure}

Figures~\ref{fig:duty_factor1} and \ref{fig:duty_factor2} show the changes in
the average achieved capacity during a day and the average increase in data
center capacity achieved as wind level grows, respectively.  We can see that a
capacity of 20\% is achieved at 0\% wind level. This indicates that
dispatchable loads are used to decrease spillages of nonwind renewables and
imports but the data-center loads are far from fully served, an unattractive
feature from a data-center owner standpoint.  Achieved data-center capacity,
however, increases rapidly to 60\% at a 5\% wind level (indicating that
stranded wind power adds value to the loads).  In addition, achieved capacity
goes up to 80\% at 50\% wind penetration level, and the trend is maintained
throughout the day. We found that the variance of the achieved capacity does
not increase for wind levels higher than 15\%. We also note that unserved loads
for flexible data centers are not penalized in Case 4 (no demand elasticity
price is provided) but high achieved data-center capacity can still be
achieved. Consequently, we conclude that stranded power provides a natural
economic incentive for flexible computing. In addition, the result indicates
that the power grid indeed benefits from serving the data center loads.

\section{Summary and Future Work}
\label{section:summary}

Our analysis shows that increased wind penetration levels lead to high levels
of spillage and uneconomic absorbed generation, which together we call stranded
power.  Significant at even moderate levels of wind penetration, these numbers
grow even higher at high levels of wind penetration (83\% at 50\% penetration).
Two scenarios added data centers in conventional ways, alone and paired with a
wind power plant.  However, a new kind of scenario, which added data centers as
a dispatchable load that the grid could turn on or off based on grid benefits,
gave surprising results.  Spillage and average power cost decreased
dramatically; as much as 44\%.  Dispatchable loads enable better utilization of
wind generation and significantly more efficient and flexible network control,
enabling dispatchable data center loads to achieve 60$\sim$80\% of full
capacity. We are planning to extend our analysis to capture more detailed
data-center scheduling models and network models of higher fidelity.



\section*{Acknowledgment}

This material is based upon work supported by the U.S. Department of Energy, 
Office of Science, under contract number DE-AC02-06CH11357 and the National 
Science Foundation under Award CNS-1405959. We gratefully acknowledge the 
computing resources provided on Blues, a high-performance computing cluster 
operated by the Laboratory Computing Resource Center at Argonne National 
Laboratory. We thank Anthony Papavasiliou for generously sharing the Western 
Electricity Coordinating Council test system data.

\ifCLASSOPTIONcaptionsoff
  \newpage
\fi





\bibliographystyle{IEEEtran}

\def\UrlBreaks{\do\/\do-\do\.\do:}
\bibliography{zccloud}

\begin{framed}
\noindent The submitted manuscript has been created by UChicago Argonne, LLC, Operator of Argonne National Laboratory (``Argonne'').  Argonne, a U.S. Department of Energy Office of Science laboratory, is operated under Contract No. DE-AC02-06CH11357.  The U.S. Government retains for itself, and others acting on its behalf, a paid-up nonexclusive, irrevocable worldwide license in said article to reproduce, prepare derivative works, distribute copies to the public, and perform publicly and display publicly, by or on behalf of the Government. 
\end{framed}

\end{document}